\documentclass[12pt]{amsart}
\usepackage{amscd,amssymb}
\usepackage[graph,frame,poly,arc]{xy}  
\usepackage[plainpages,backref,urlcolor=blue]{hyperref}

\theoremstyle{plain}
\newtheorem{thm}[subsection]{Theorem}
\newtheorem{lem}[subsection]{Lemma}
\newtheorem{prop}[subsection]{Proposition}
\newtheorem{cor}[subsection]{Corollary}

\theoremstyle{definition}
\newtheorem{rk}[subsection]{Remark}
\newtheorem{definition}[subsection]{Definition}
\newtheorem{ex}[subsection]{Example}

\newtheorem{question}[subsection]{Question}

\numberwithin{equation}{section}
\newcommand{\OO}{{\mathcal O}}

\newcommand{\F}{{\mathcal F}}

\newcommand{\al}{{\alpha}}
\newcommand{\be}{{\beta}}

\newcommand{\Z}{\mathbb{Z}}

\newcommand{\C}{\mathbb{C}}
\newcommand{\PP}{\mathbb{P}}

\begin{document}
%\date{June 4, 2009}

\title [Ramblings on the freeness of affine  hypersurfaces]
{Ramblings on the freeness of affine  hypersurfaces}

\author[Alexandru Dimca]{Alexandru Dimca$^{1}$}
\address{Universit\'e C\^ ote d'Azur, CNRS, LJAD and INRIA, France and Simion Stoilow Institute of Mathematics,
P.O. Box 1-764, RO-014700 Bucharest, Romania}
\email{dimca@unice.fr}

\author[Gabriel Sticlaru]{Gabriel Sticlaru}
\address{Faculty of Mathematics and Informatics,
Ovidius University
Bd. Mamaia 124, 900527 Constanta,
Romania}
\email{gabriel.sticlaru@gmail.com}

\thanks{$^1$ This work has been partially supported by the Romanian Ministry of Research and Innovation, CNCS - UEFISCDI, grant PN-III-P4-ID-PCE-2020-0029, within PNCDI III.
}



\begin{abstract}  In this note we look at the freeness for complex affine hypersurfaces. If $X \subset \mathbb{C}^n$ is such a hypersurface, and $D$ denotes the associated projective hypersurface, obtained by taking the closure of $X$ in $\mathbb{P}^n$, then we  relate first the Jacobian syzygies of $D$ and those of $X$. Then we introduce two types of freeness for an affine hypersurface $X$, and prove various relations between them and the freeness of the projective hypersurface $D$. We write down a proof of the folklore result saying that an affine hypersurface is free if and only if all of its singularities are free, in the sense of K. Saito's definition in the local setting. 
In particular, smooth affine hypersurfaces and affine plane curves are always free.
Some other results, involving global Tjurina numbers and minimal degrees of non trivial syzygies are also explored.
\end{abstract}
 
\maketitle

%\tableofcontents

\section{Introduction} 

The notion of free hypersurface, or free divisor, has a long and fascinating history. This concept was considered until now mostly in two settings, namely
\begin{enumerate}

\item[{(L)}] The local setting, where one looks at germs of hypersurfaces, as in the seminal paper by K. Saito \cite{KS}, see Definition \ref{defFreeLoc} below. In particular, any plane curve singularity is free.

\item[{(P)}]  The projective setting, where one considers hypersurfaces in a complex projective space $\PP^n$, as for instance in \cite{B+,Dmax,DStFD,DStFS,Sim2,ST,To}. In particular, the free projective hypersurfaces satisfy a lot of restrictions, for instance they must be rather singular, namely the singular locus must have codimension one, see \cite[Theorem 2.8]{DStFD} for curves and
\cite[Corollary 4.4]{DStFS}  for surfaces.
\end{enumerate}

Note that the  central hyperplane arrangements in some affine space $\C^m$ can be considered from both view-points, and Terao's Conjecture, see for instance \cite{DHA}, is a beautiful open problem in this area.
More generally, any hypersurface in $\C^n$ defined by a homogeneous polynomial can be considered from both view-points, and freeness in one setting is the same as freeness in the other.

\medskip

In this paper we discuss the freeness for a reduced hypersurface $X:g=0$ in $\C^n$, where $g \in R=\C[x_1,\ldots,x_n]$, as well as a number of related results. The main object of study is the  $R$-module $AR(g)$ of syzygies among the partial derivatives
$g_1, \ldots,g_n$ of $g$ with respect to $x_1, \ldots, x_n$ and $g$ itself. This $R$-module  $AR(g)$ is isomorphic to the $R$-module of  $Der(g)$ of derivations of $R$ preserving the principal ideal $(g) \subset R$, see Definition \ref{def2} and Lemma \ref{ldef}.
However,  $AR(g)$ seems to contained more information than $Der(g)$, for instance it allows us to define a {\it syzygy-degree function}
$sdeg:AR(g) \to \Z$, see Definition  \ref{defdeg}, and the numerical invariant $mdr(g)$, the minimal degree of a non-trivial syzygy in $AR(g)$. We denote by $Der_0(g)$ the $R$-submodule in $Der(g)$ consisting of those derivations killing the equation $g$.

If $D:f=0$ is the usual projectivization of the hypersurface $X$, obtained by taking the closure of $X$ in the complex projective space $\PP^n$, then in Theorem \ref{thm1} we show that there is a natural epimorphism $\phi: AR(f) \to AR(g)$, where $AR(f)$ is the module of syzygies among the partial derivatives
$f_0, \ldots,f_n$ of $f$. This result is refined in Corollary \ref{corthm1}, where we show that there are isomorphisms
$\phi_s: AR(f)_s \to AR(g)_{\leq s}$ for all integers $s$.

For affine hypersurfaces, we introduce two notions of freenees, the $A$-freeness and the $A_0$-freeness, related to the freeness of the $R$-module $Der(g)$, and respectively $Der_0(g)$, see Definition \ref{defFree}. The $A$-freeness is close to the usual freeness in the setting (L) and the $A_0$-freeness is close to the usual freeness in the setting (P). In Proposition \ref{propF0} and Corollary \ref{corF0} we describe some situations when the affine hypersurface $X$ is $A$-free if and only if $X$ is $A_0$-free. The following question is open, though we believe the answer should be affirmative.

\begin{question}
\label{q1} 
Does there exist an affine hypersurface $X \subset \C^n$, with $n \geq 3$, such that $X$ is $A_0$-free, but not $A$-free ?
\end{question}

There is a characterization of bases of the free $R$-module $Der(g)$, the so called {\it Saito's Criterion}, similar to the known results in the the settings (L) and (P), see Proposition \ref{propA}. Using this result, we show that the {\it Milnor fibers of weighted homogeneous polynomials} are always $A$-free, and sometimes also $A_0$-free, see Proposition
\ref{propKS}, and construct explicit bases for the free $R$-module $Der(g)$ in these cases.

Another class of free affine hypersurfaces comes from the following fact: if the projective hypersurface $D$ is free, then the affine
part $X$ is $A$-free, see Corollary \ref{corF1}.

We write down a proof of the folklore result saying that an affine hypersurface is free if and only if all of its singularities are free, in the sense of K. Saito's definition in the local setting, see Theorem \ref{thmA}. A rather a surprizing consequence of this result is that any smooth affine hypersurface is $A$-free, while it was known that a free projective hypersurface must have a codimension one singular locus.

In particular, any affine plane curve $X:g=0$ is both $A$-free and $A_0$-free, see Corollary \ref{corA}. Moreover, in this case the ideal $A(g)=J_g:(g)$, where  $J_g$ denotes the Jacobian ideal of $g$,  is either the whole ring $R$, or it is minimally generated by exactly two elements. 

In Corollary \ref{corthm2} we give an upper bound for the total Tjurina number of an affine plane curve $X$ in terms of the minimal syzygy-degree $mdr(g)$ of a syzygy in $AR(g)$, a result obtained in the projective setting (P) by 
A. du Plessis and C. T. C. Wall, see \cite{dPW} and then reproved by several authors, see for instance \cite{Dmax,E}.
A similar, but somewhat weaker result holds for any hypersurface $X$ such that its closure $D$ has only isolated singularities, see Corollary \ref{corthmC}.
Another relation between the total Tjurina number of the projective hypersurface $D$ with only isolated singularities and $mdr(g)$ is given in Proposition \ref{propB}.

\medskip

We would like to thank Laurent Bus\' e and Masahiko Yoshinaga for useful discussions concerning this project.

\section{Affine vs. projective syzygies and their degrees} 

Let $R=\C[x_1,\ldots,x_n]$ be the polynomial ring in $n$ variables  with complex coefficients, and let $X:g=0$ be a reduced hypersurface of degree $d$ in the complex affine space $\C^n$.

Let $J_g=(g_1,\ldots,g_n)\subset R$ denote the corresponding {\it Jacobian ideal} where $g_j= \partial_jg$. Here and in the sequel
$$\partial_j= \frac{\partial }{\partial x_j} \text{ for all } j=1, \ldots, n.$$

One can introduce the {\it Tjurina ideal} of $g$ by setting
$$T_g=J_g+(g)$$
where $(g)$ denotes the principal ideal generated by $g$ in $R$, and then
the singular set $X_{sing}$ of $X$ is precisely the zero set of this ideal $T_g$.

 \begin{rk}
\label{rkC1} When $n=2$, $X$ being reduced is equivalent to  $X$ having only isolated singularities, and this is in turn equivalent to the fact that the partial derivatives $g_1$ and $g_2$ form a regular sequence in $R$.
Let $\tau(X)$ denote the {\it total Tjurina number} of the affine curve $X$, which is by definition the sum of all local Tjurina numbers $\tau(X,p)$ for $p \in X_{sing}$. Then one clearly has
\begin{equation}
\label{tau} 
 \tau(X)= \dim \frac{R}{T_g}.
\end{equation}
 \end{rk}
 
Note that in the general case there is an exact sequence
\begin{equation}
\label{es1} 
0 \to AR(g) \to R^{n+1} \to T_g \to 0,
\end{equation}
where the morphism
$\psi:R^{n+1}\to T_g$ is given by 
\begin{equation}
\label{ideal1.1} 
\psi(\al_0, \al_1, \ldots, \al_n)= \al_0 g +\al_1 g_1 +\ldots + \al_n g_n,
\end{equation}
$AR(g)=\ker \psi$ and the morphism $AR(g) \to R^{n+1}$ is the inclusion.
In addition, to the hypersurface $X:g=0$ we can associated the ideal
\begin{equation}
\label{ideal1} 
 A(g)= J_g:(g)=\{\be \in R \ : \  \be g \in J_g \} \subset R.
\end{equation}
 There is an exact sequence of $R$-modules
\begin{equation}
\label{ideal1.2} 
0 \to AR_0(g) \to  AR(g) \to A(g) \to 0, 
\end{equation}
where the second morphism $\pi: AR(g) \to A(g)$ is the projection 
$$(\al_0, \al_1, \ldots, \al_n) \mapsto \al_0,$$
$AR_0(g)= \ker \pi$, and $AR_0(g) \to  AR(g)$ is the inclusion.
\begin{ex}
\label{exC1} 
When $n=2$, the $R$-module $AR_0(g)$ is free of rank 1, generated by
$(0,g_2,-g_1)$. This follows from the fact that in this case $g_1$ and $g_2$ have no common factors, see Remark \ref{rkC1}.
\end{ex}
Let $Der(R)=\{\partial =\al_1 \partial_1 + \ldots + \al_n \partial_n \ : \ \al_j  \in R \}$ be the free $R$-module of $\C$-derivations of the polynomial ring $R$.
\begin{definition}
\label{def2} 
The $R$-module $Der(g)$ of derivations of $R$ preserving the principal ideal $(g) \subset R$ is by definition
$$Der(g)=\{\partial \in Der(R) \ : \ \partial g \in (g) \}.$$
Moreover, the $R$-module $Der_0(g)$ of derivations of $R$ killing the equation $g$  is by definition
$$Der_0(g)=\{\partial \in Der(R) \ : \ \partial g =0 \}.$$
\end{definition}
There are similar definitions in the local case, where the polynomial ring $R$ is replaced by the local ring $\OO_n$ of germs of analytic functions at the origin of $\C^n$. Let $Der(\OO_n)=\{\partial =\al_1 \partial_1 + \ldots + \al_n \partial_n \ : \ \al_j  \in \OO_n \}$ be the free $\OO_n$-module of $\C$-derivations of the ring $\OO_n$.
\begin{definition}
\label{def2loc} 
Consider a hypersurface singularity $(X,0)$ at the origin of $\C^n$, given by an equation $h=0$, with $h \in \OO_n$.
The $\OO_n$-module $Der(h)$ of derivations of $\OO_n$ preserving the principal ideal $(h) \subset \OO_n$ is by definition
$$Der(h)=\{\partial \in Der(\OO_n) \ : \ \partial h \in (h) \}.$$
Moreover, the $\OO_n$-module $Der_0(h)$ of derivations of $\OO_n$ killing the equation $h$  is by definition
$$Der_0(h)=\{\partial \in Der(\OO_n) \ : \ \partial h =0 \}.$$
\end{definition}
However, note that the module $Der(h)$ does not depend on the chosen equation for the singularity $(X,0)$, but the module $Der_0(h)$ does depend, so the interest of the latter is rather limited in the local setting.

Coming back to the global, polynomial case, one has the following.
\begin{lem}
\label{ldef}
The $R$-modules $Der(g)$ and $AR(g)$ are isomorphic.
Similarly, the $R$-modules $Der_0(g)$ and $AR_0(g)$ are isomorphic.
\end{lem}
\proof
It is clear that the obvious projection, namely 
$$AR(g) \to Der(g), \ (\al_0, \al_1, \ldots, \al_n) \mapsto  \partial =\al_1 \partial_1 + \ldots + \al_n \partial_n, $$
is an isomorphism. The claim for the $R$-modules $Der_0(g)$ and $AR_0(g)$ is even more obvious.
\endproof

\bigskip

Let $S=\C[x_0,x_1,\ldots,x_n]$ be the graded polynomial ring in $n+1$ variables  with complex coefficients, and let $D:f=0$ be a reduced hypersurface of degree $d$ in the complex projective plane $\PP^2$. 
We denote by $J_f$ the {\it Jacobian ideal} of $f$, i.e. the homogeneous ideal in $S$ spanned by the partial derivatives $f_0,f_1,\ldots, f_n$ of $f$, where $f_0=\partial_0f$ and $\partial_0$ are defined in the obvious way.
Consider the graded $S$-module of {\it Jacobian syzygies} of $f$, namely
$$AR(f)=\{(a_0, a_1, \ldots, a_n) \in S^{n+1} \ : \  a_0f_0+a_1f_1+\ldots +a_nf_n=0\}.$$

If one starts with an affine hypersurface $X:g=0$ as above, one can consider its closure $D:f=0$ in $\PP^n$ under the inclusion $\C^n \to \PP^n$ given by 
$$(x_1,\ldots,x_n) \mapsto (1:x_1:\ldots: x_n).$$
 Algebraically, this means that
$$f(x_0,x_1,\ldots,x_n)=x_0^dg\left(\frac{x_1}{x_0}, \ldots, \frac{x_n}{x_0}\right) \text{ and } g(x_1,\ldots,x_n)=f(1,x_1,\ldots,x_n).$$
Note that the homogeneous polynomial $f$ obtained in this way is not divisible by $x_0$ or, in other words, the hyperplane at infinity $H_0:x_0=0$ is not an irreducible component for $D$. In the sequel we consider only pair of hypersurfaces $(X,D)$ obtained by this natural construction.

Now we relate the $S$-module $AR(f)$ to the $R$-module $AR(g)$.
Note that any $R$-module can be regarded as an $S$-module using the
ring morphism
$$\theta: S \to R, \  \  \theta(u )=u(1,x_1,\ldots,x_n),$$
for any $u\in S$. Note that the $S$-modules obtained using $\theta$ are not graded $S$-modules in general. There is an obvious  $\C$-linear isomorphism
$$\eta_e: R_{\leq e}  \to S_e$$
for any integer $e \geq 0$, given by
$$\eta_e(v)= x_0^ev\left(\frac{x_1}{x_0}, \ldots, \frac{x_n}{x_0}\right).$$
 Here $R_{\leq e}= \{ v \in R \ : \ \deg v \leq e\}$ and $S_e$ is the homogeneous component of degree $e$ of the polynomial ring $S$.
 The following result is obvious.
\begin{lem}
\label{lem2}
If $v \in R_{\leq e}$ and $v' \in R_{\leq e'}$, then
$$\eta_{e+e'}(vv')=\eta_e(v) \eta_{e'}(v').$$
\end{lem}

\begin{thm}
\label{thm1}
For any pair of hypersurfaces $X:g=0$, $D:f=0$ as above, one has a surjective morphism of $S$-modules 
$\phi: AR(f) \to AR(g)$, given by
$$ \phi (a_0, a_1, \ldots, a_n) =( d\cdot \theta(a_0), \theta (a_1)-x_1\theta(a_0),\ldots, \theta(a_n)-x_n\theta(a_0)), $$
 for any  $(a_0, a_1, \ldots, a_n) \in AR(f).$

\end{thm}

\proof
First we show that the morphism $\phi$ is well-defined, that is that one has $\phi (a_0, a_1, \ldots, a_n) \in AR(g)$.
The equality $g(x_1,\ldots,x_n)=f(1,x_1,\ldots,x_n)$ implies the following
\begin{equation}
\label{eq1} 
 g_j=\theta(f_j) \text{ for any } j=1, \ldots,n.
\end{equation} 
The Euler formula for the homogeneous polynomial $f$, namely
$$x_0f_0+x_1f_1+\ldots + x_nf_n=d\cdot f,$$
implies the equality
\begin{equation}
\label{eq2} 
 \theta(f_0) =d \cdot g-x_1g_1-\ldots -x_ng_n.
\end{equation} 
Using these formulas, we see that applying $\theta$ to the equality
$$ a_0f_0+a_1f_1+\ldots +a_nf_n=0$$
we get
\begin{equation}
\label{eq3} 
d\cdot \theta(a_0)g+
(\theta (a_1)-x_1\theta(a_0))g_1+\ldots + (\theta (a_n)-x_n\theta(a_0))g_n =0.
\end{equation} 
Therefore the morphism $\phi$ is well-defined. Now we show that the morphism $\phi$ is surjective. Let  $\rho'=(\al_0, \al_1, \ldots, \al_n) \in R^{n+1}$ be a vector such that
\begin{equation}
\label{eq3.1} 
 \al_0 g +\al_1 g_1 +\ldots + \al_n g_n=0
\end{equation} 
Comparing the equations \eqref{eq3} and \eqref{eq3.1}, we see that it is enough to find a vector 
$\rho=(a_0, a_1, \ldots, a_n) \in AR(f)$ such that
\begin{equation}
\label{eq3.2} 
\alpha_0=   d \cdot \theta (a_0)
 \text{ and }   \al_j= \theta (a_j)-x_j\theta(a_0), \text{ for all } 1 \leq j \leq n.
\end{equation} 
This system of equations has a unique solution
\begin{equation}
\label{eq3.3} 
\theta (a_0)=\frac{\al_0}{d}  \text{ and }
 \theta (a_j)= \al_j+\frac{x_j\al_0}{d},   \text{ for all } 1 \leq j \leq n.
\end{equation} 
In other words, the vector  $(\theta (a_0),\cdots, \theta (a_n))$, if it exists, is uniquely determined by the vector $\rho'$.
This solution satisfies the equation
\begin{equation}
\label{eq3.4} 
 \theta (a_0) (d \cdot g-x_1g_1-\ldots- x_ng_n)+
 \theta (a_1)g_1+ \ldots +\theta (a_n)g_n =0.
\end{equation} 
We set
\begin{equation}
\label{eq4}
s=\max\{ \deg \al_j \ : \ 1 \leq j \leq n \}.
\end{equation} 
Note that equation \eqref{eq3.1} implies that 
$$\deg( \al_0 g) \leq \max\{ \deg(\al_j g_j) \}\leq s+d-1,$$
and hence $\deg \al_0 \leq s-1$. It follows that
\begin{equation}
\label{eq4.1}
\max\{ \deg  \al_0,  \deg (  d \cdot \al_j +x_j \al_0 ) \text{ for } 1 \leq j \leq n  \} \leq s.
\end{equation} 

Using \eqref{eq1} and \eqref{eq2}, we see that
\begin{equation}
\label{eq5}
f_j=\eta_{d-1}(g_j) \text{ for } 1 \leq j \leq n ,   \text{ and } f_0= \eta_{d-1}(d \cdot g-x_1g_1-\ldots- x_ng_n).
\end{equation}
To end the proof we apply Lemma \ref{lem2} with $e=s$ and $e'=d-1$ to the equality \eqref{eq3.4}. In this way we get a vector
$$\rho=(a_0, a_1, \ldots, a_n) \in AR(f)_s \text{ such that  } \phi(\rho)=\rho',$$
where $a_j=\eta_s( \al_j+ \frac{x_j\al_0}{d})$ for $ 1 \leq j \leq n $ and 
$ a_0= \eta_s( \frac{\al_0}{d}). $
If the inequality in \eqref{eq4.1} is strict, then the syzygy $\rho$ can be simplified by a power of $x_0$, yielding a new syzygy $\rho_1$ of degree $<s$ such that $\phi(\rho_1)=\rho'$.
\endproof
Motivated by the above proof we introduce the following.
\begin{definition}
\label{defdeg}
We define the {\it syzygy-degree} $sdeg(\rho')$ for any non-zero element $$\rho'=(\al_0, \al_1, \ldots, \al_n) \in AR(g)$$
by the formula
$$sdeg(\rho') =  \max\{ \deg  \al_0,  \deg (  d \cdot \al_j +x_j \al_0 ) \text{ for } 1 \leq j \leq n  \} .$$
\end{definition}
Note that one has 
\begin{equation}
\label{eq6}
\deg  \al_0 \leq sdeg(\rho') \leq  \max \{\deg \al_0 +1,  \deg \al_j \text{ for } 1 \leq j \leq n \},
\end{equation}
and both these inequalities can be strict, see Remark \ref{rkst}.
Note also the following obvious but useful property: for any non-zero polynomial $u \in R$ and  any non-zero element $\rho' \in AR(g)$, one has
\begin{equation}
\label{eq6.1}
sdeg(u \rho')= \deg u + sdeg(\rho').
\end{equation}

With this notation, we have the following result.

\begin{cor}
\label{corthm1} 
With the above notation, for any pair of hypersurfaces $X:g=0$, $D:f=0$ as above and for any integer $s$, there is an isomorphism of $\C$-vector spaces
$$\phi_s : AR(f)_{s} \to AR(g)_{\leq s},$$
where $AR(g)_{\leq s} =\{(\rho'  \in AR(g) \ : \ sdeg(\rho')  \leq s\}$. 
In particular, if we define
\begin{equation}
\label{eq7}
mdr(f)= \min\{r \ : \ AR(f)_r \ne 0\} \text{ and } mdr(g)= \min\{r \ : \ AR(g)_{\leq r} \ne 0\},
\end{equation}
then $$
mdr(g)=mdr(f)\leq d-1.$$

\end{cor}
\proof  The proof of Theorem \ref{thm1} yields the surjectivity of $\phi_s$.
To prove the injectivity of $\phi_s$, let $\rho \in AR(f)_s$ be a  syzygy such that
$\phi(\rho)=0$. This yields $\theta(a_j)=0$ for all $0 \leq j \leq n$, where $\rho=(a_0, a_1, \ldots, a_n)$. This implies $\rho$, since the restriction of $\theta$ to each homogeneous component $S_k$ is clearly injective. 
The inequality $mdr(f) \leq d-1$ comes from the existence of the Koszul syzygies $\kappa_{ij}$, see below.
\endproof
The following result is a direct consequence of the last claim in Corollary \ref{corthm1}.
\begin{cor}
\label{cor10} 
Let $D:f=0$ be a  hypersurface of degree $d>1$ in $\PP^n$. Consider all the affine parts of it, namely the affine hypersurfaces
$$ X_H=D \setminus H \subset \PP^n \setminus H \simeq \C^n,$$
for all hyperplanes $H \subset \PP^n$, $H$ not an irreducible component of $D$. If $g_H=0$ is an equation for the hypersurface $X_H$, then one has $mdr(g_H)=mdr(f)$ for all these hyperplanes $H$.
\end{cor}

Note that among the Jacobian syzygies in $AR(f)$ there are the following {\it Koszul syzygies} $\kappa_{ij}$, for all $0 \leq  i <j \leq n$, where
\begin{equation}
\label{eq7.05}
\kappa_{ij}=f_ie_j -f_je_i,
\end{equation}
where the vector $e_k \in \C^{n+1}$ has 1 on the $k$-th coordinate and zero everywhere else, for $k=0,\ldots,n$.
%They satisfy one second order syzygy, namely
%\begin{equation}
%\label{eq7.1}
%f_x\kappa_x+f_y\kappa_y+f_z\kappa_z=0.
%\end{equation}

Define new syzygies in $AR(g)$ by 
$\kappa'_{ij}=\phi(\kappa_{ij}),$ using the morphism $\phi$ from Theorem  \ref{thm1}.
Note that for $0<i<j\leq n$, one has
$\kappa'_{ij}=g_ie_j -g_je_i$, while for $i=0$ and $1 \leq j \leq n$ one has
$$\kappa'_{0j}=-dg_je_0+ g_j\sum_{k=1,n; k \ne j}x_ke_k+(d \cdot g -\sum_{k=1,n; k \ne j}x_kg_k)e_j.$$

\begin{rk}
\label{rkst}
Note that the syzygies $\kappa'_{ij}$ give examples where the inequalities in \eqref{eq6} are strict.
\end{rk}

Let $KR(f)$ be the $S$-submodule in $AR(f)$ generated by the Koszul syzygies $\kappa_{ij}$, and $KR(g)$ be the $R$-submodule in $AR(g)$ generated by the syzygies $\kappa'_{ij}$.
Consider the quotient modules
\begin{equation}
\label{eq9}
ER(f)=AR(f)/KR(f), \ ER(g)=AR(g)/KR(g)  \text{ and } E(g)=A(g)/J_g.
\end{equation}

Note that the polynomial ring $R$ has a natural increasing filtration
 $$R_{\leq s}=\{ h \in R \ : \  \deg h \leq s \}$$
 and this filtration gives rise to an increasing filtration on the ideal $A(g) \subset R$,  namely
$A(g)_{\leq s}= A(g) \cap R_{\leq s}$, and also
 on the quotient  ideal $E(g)$, namely
 $$E(g)_{\leq s}=\{ [h] \in E(g) \ : \   h \in A(g)_{\leq s}\}.$$
 Similarly, the $sdeg$-filtration $AR(g)_{\leq s}$ on the module $AR(g)$ induces an increasing filtration on the quotient $ER(g)$, namely
$$ER(g)_{\leq s}=\{ [\rho'] \in ER(g) \ : \   \rho' \in AR(g)_{\leq s}\}.$$ 
\begin{cor}
\label{cor1.1} 
With the above notation, for any pair of curves $X:g=0$, $C:f=0$ as above, there are  surjective morphisms
$\phi: ER(f) \to ER(g)$ and  $ \pi: ER(g) \to E(g).$
In addition, one has a surjective morphism
$$ \phi_s: ER(f)_s \to ER(g)_{\leq s},$$
for any integer $s$.
\end{cor} 
\proof
First we explain the notation. The morphism $\phi: ER(f) \to ER(g)$ is induced by the morphism $\phi: AR(f) \to AR(g)$ and the surjectivity 
of $\phi: ER(f) \to ER(g)$ and of $\phi_s: ER(f)_s \to ER(g)_{\leq s} $
follows from Theorem \ref{thm1} and Corollary \ref{corthm1} and the definition of the quotient modules $ER(f)$ and $ ER(g)$. 

The morphism $\pi: ER(g) \to E(g)$ is induced by the projection on the third factor $\pi: AR(g) \to A(g)$, and
 the surjectivity follows the definition of the quotient modules $ER(g)$ and $ E(g)$.

 \endproof

 \begin{rk}
\label{rk1.0}
Concerning the morphisms in Corollary \ref{cor1.1} there are several useful points to make.

\begin{enumerate}

\item For the case of curves, namely when $n=2$, the morphism $ \pi: ER(g) \to E(g)$ is an isomorphism. To show that $ \pi: ER(g) \to E(g)$ is injective, assume that
 $\pi ([\al_0,\al_1,\al_2])=0$ in $E(g)$. This means that $\al_0 \in J_g$, hence one can write $\al_0=u_1g_1+u_2g_2$ for some $u_1,u_2 \in R$.
 This implies that 
 $$(\al_0,\al_1,\al_2)+\frac{u_1}{d}k'_{01}+\frac{u_2}{d}k'_{02} \in \ker \{\pi: AR(g) \to A(g)\}=AR_0(g),$$
 and then use the description of $AR_0(g)$ given in Example \ref{exC1}.
 It follows that $[\al_0,\al_1,\al_2]=0$ in $ER(g)$, which proves the claimed injectivity.

\item  It is not true that the restrictions
 $\pi : ER(g)_{\leq s} \to E(g)_{\leq s}$
 are surjective.
 Indeed, if $g=x_1^2+x_2^3$, then $g \in J_g$ and hence $A(g)=R$ and $E(g)_0 \ne 0$. On the other hand, $ER(g)_{\leq 0} \ne 0$ would imply by Corollary \ref{cor1.1} that $ER(f)_0 \ne 0$. But a homogeneous polynomial $f$ satisfies $AR(f)_0 \ne 0$ if and only if after a coordinate change $f$ does not depend on one variable, and hence the zero set $C:f=0$ is a union of lines passing through one point. This shows that $ER(g)_{\leq 0} = 0$ when $X:g=0$ is not a family of lines.

\item It is not true that the restrictions $\phi_s: ER(f)_s \to ER(g)_{\leq s} $ are injective. If $f=x_0^2x_2+x_1^3$ and hence $g=x_2+x_1^3$, one has
that the syzygy $\rho=(x_0,0, -2x_2) \in AR(f)_1$ is not in $KR(f)_1=0$.
However, $\rho'=\phi_1(\rho)=(3,-x_1,-3x_2) \in KR(g)_1$.
More precisely, one has
$$\phi_2(-\kappa_{02})=-\kappa'_{02}=(3,-x_1,-3x_2)$$
and
$sdeg(3,-x_1,-3x_2)=1$, which proves our claim.

\end{enumerate}
 
  \end{rk}

\begin{rk}
\label{rk1.1}
The interest of the last claim in Corollary \ref{cor1.1} comes from the fact that we have some strong vanishing results for $ER(f)_s$, which translate in having only high degree representatives for the classes in $ER(g)$ in many cases. 
For instance, if $D$ has only simple nodes $A_1$ as singularities, then it is known that $ER(f)_s =0$ for $s \leq nd/2-n-1$, see \cite{DMS,DS,DStEdin}.
When the hypersurface $D$ has a single node, then it is known that $ER(f)_s=0$ for
$s <n(d-2)$, see \cite[ Definition 1.1, Equation 1.3 and Example 4.3]{DStEdin}, hence the generators of $ER(g)$ have $sdeg$ at least $n(d-2)$.
\end{rk}
\begin{rk}
\label{rk1.2}
Both for the affine hypersurface $X:g=0$ and for the projective hypersurface $D:f=0$, the defining equation is uniquely determined up to a non-zero constant factor. It follows that the objects defined above, namely $J_g,T_g, AR(g), AR_0(g), KR(g), ER(g), A(g)$, $Der(g), Der(g)_0$, $mdr(g)$ can be denoted when convenient simply by $J_X,T_X, AR(X), AR_0(X), KR(X), ER(X)$, $A(X), Der(X),$
$ Der_0(X), mdr(X)$, since they depend only on $X$ and not on the chosen equation $g=0$ for $X$. Similarly, $J_f, AR(f), KR(f), ER(f), mdr(f)$ can be denoted simply by 
$J_D, AR(D), KR(D), ER(D), mdr(D)$.
\end{rk}

\section{On the freeness of  affine  hypersurfaces} 

First we recall the following definition, given by K. Saito in \cite{KS}. 

\begin{definition}
\label{defFreeLoc} 
A  hypersurface singularity $(X,0):h=0$ at the origin of $\C^n$ is said to be free if the $\OO_n$-module of derivations $Der(h)$ introduced in Definition \ref{def2loc} is free.
\end{definition}
Motivated by this local notion, we introduce the following two global, algebraic notions.
\begin{definition}
\label{defFree} 
The reduced affine hypersurface $X \subset \C^n$ is $A$-free (resp. $A_0$-free)
if the $R$-module $Der(X)$ of derivations preserving $X$ (resp.
the $R$-module $Der_0(X)$ of derivations killing an equation of  $X$) is free.
\end{definition}
Using the exact sequences \eqref{es1} and  \eqref{ideal1.2}, and Lemma \ref{ldef}, one sees that the rank of the $R$-module $Der(X)$ is $n$, while the rank of 
the $R$-module $Der_0(X)$ is $n-1$. Hence the $R$-module $Der(X)$ (resp. $Der_0(X)$)
is free if and only if it can be generated by $n$ (resp. $n-1$) derivations.

The $A$-freeness defined above is similar to the usual definition in the local setting (L) given by K. Saito in \cite{KS}. 
In the projective setting (P), that is when $g$ is itself a homogeneous polynomial, the $A$-freeness and the $A_0$-freeness are equivalent,  see for instance \cite[Definition 8.1]{DHA}. In fact, we have the following more general result.
\begin{prop}
\label{propF0}
Let $X:g=0$ be an affine hypersurface such that there is a derivation
$$\epsilon=\al_1\partial_1+ \ldots + \al_n \partial_n \in Der(X)$$
such that $\epsilon g= c \cdot g$, where $c \in \C^*$.
Then there is a direct sum decomposition
$$Der(X)=Der_0(X) \oplus R\cdot \epsilon,$$
and $X$ is $A$-free if and only if $X$ is $A_0$-free.
\end{prop}
Note that the hypothesis in Proposition \ref{propF0} is equivalent to the equality $A(X)=R$.
\proof
It is clear that $\epsilon \in Der(g)=Der(X)$, and that for any $\delta \in Der(g)$, satisfying $\delta g=\al \cdot g$ for some $\al \in R$, one has
$$\delta - \al  \cdot \epsilon \in Der_0(g)=Der_0(X).$$
This  yields the claimed direct sum decomposition. Then, if $Der_0(x)$ is a free $R$-module, then this decomposition implies that  $Der(X)$
is also free. Conversely, if $Der(X)$ is free, then $Der_0(X)$ is a projective $R$-module, being a direct summand in a free module.
We get the result using  Quillen-Suslin Theorem, saying that any finite type projective $R$-module is free, see \cite{Q,Su}.
\endproof

\begin{cor}
\label{corF0}
An affine hypersurface $X:g=0$ is $A$-free if and only if it is $A_0$-free,
if one of the following cases occurs.
\begin{enumerate}

\item The polynomial function $g: \C^n \to \C$ has no singularities, i.e. it is a submersion. Equivalently, $J_g=R$.

\item The polynomial $g$ is weighted homogeneous of type $(w_1,\ldots, w_n;N)$, where the weights $w_j$ are positive or negative integers, and the weighted degree $N$ is not zero.

\end{enumerate}
\end{cor}
In the case (1) above, one has in fact $J_g=A_g=A(X)=R$, and such polynomials $g$ can be rather complicated, in particular far away from weighted homogeneous polynomials, see \cite{ACL}.
Note that for the claim (2), the Euler derivation is
\begin{equation}
\label{Euler}
\epsilon=\frac{w_1x_1}{N}\partial_1+ \ldots + \frac{w_nx_n}{N} \partial_n .
\end{equation}

For general affine hypersurfaces, the relation between 
$A$-freeness and $A_0$-freeness is not clear, as we see in Proposition  \ref{propKS} and in Theorem \ref{thmA} below.

In order to find a basis for the free $R$-module $AR(g)$, the following result might be useful. Note that its local analog is well known, see \cite{KS}, as well as the projective analog, see \cite[Theorem 8.1]{DHA}.

\begin{prop}(Saito's Criterion)
\label{propA} 
Let $\delta_i=\al_{i1}\partial_1+\ldots + \al_{in} \partial_n$, for $i=1,\ldots,n$, be $n$ derivations in $Der(X)$, for the affine hypersurface $X:g=0$. Then
$\delta_1, \ldots, \delta_n$ is a basis for the free $R$-module $Der(X)$ if and only if 
$$\det M(\delta_1, \ldots, \delta_n)= \det (\al_{ij})=c \cdot g,$$
 where $c \in \C^*$. 

\end{prop}
\proof Start with any set of $n$ derivations in $Der(X)$, say $\delta_1, \ldots, \delta_n$.
If $p \in X$ is a smooth point, then one has
$$\sum_{j=1}^n\al_{ij}(p)g_j(p)=0 \text{ for } i=1,\ldots,n.$$
Since the hypersurface $X$ is reduced, the smooth points are dense on $X$, and the above  equalities imply that $\det M(\delta_1, \ldots, \delta_n)$ vanishes on $X$. Therefore 
$$\det M(\delta_1, \ldots, \delta_n)=hg,$$
for some polynomial $h \in R$. 

Assume now that $\delta_1, \ldots, \delta_n$ is a basis for the free $R$-module $Der(X)$. Then they are $R$-linearly independent, and hence $h$ is not the zero polynomial. Consider the following derivations $\be_j \in Der(g)$:
$\be_1=g\partial _1$ and then the Koszul derivations 
$\be_j=\kappa_{1j}=g_1e_j-g_je_1$ for $j=2, \ldots,n$. We derive the existence of an $n \times n$ matrix $N$ such that
$$M(\be_1, \ldots, \be_n)=N M(\delta_1, \ldots, \delta_n).$$
By taking the determinants, we get
$$gg_1^{n-1}=(\det N) h g.$$
This implies that $h$ divides $g_1^{n-1}$. By using the Koszul relations obtained by fixing a coordinate $k \ne 1$, we get in a similar way that 
$h$ divides $g_k^{n-1}$ for any $k=1, \ldots,n$. 
If we consider now the derivations $\be_j'=g\partial_j$, we get in a similar way that $h$ divides $g^{n-1}$. This gives a contradiction, since
$X:g=0$ is a reduced hypersurface.

Conversely, assume $h\in \C^*$ and let $\be$ be a derivation in $Der(X)$. Then using Cramer's rule, we see that we can find polynomial $u_i \in R$ such that
$$g \be=u_1 \delta_1 + \ldots + u_n \delta_n.$$
Let $M_i$ be the matrix obtained from $M=M(\delta_1, \ldots, \delta_n)$ by replacing the $i$-th line by the coordinates of $\be$. Then one has as above 
$\det M_i=v_i g$, for some polynomials $v_i \in R$.
Now $g \det M_i=v_i g^2$ can be regarded as the determinant of the matrix $M_i'$ obtained from $M$ by replacing the $i$-th line by the coordinates of $g\be$. Using the formula above, we see that
$\det M_i'=u_i\det M=u_ihg.$ This implies that all the coefficients $u_i$ are divisible by $g$, and hence the derivations $\delta_1, \ldots, \delta_n$ generate the $R$-module $Der(X)$.
\endproof
Using Saito's Criterion we can construct a lot of $A$-free hypersurfaces, which are sometimes also $A_0$-free. Before doing this, we recall  that the hypersurface $D:f=0$ in $\PP^n$ is  {\it free} if the graded $S$-module $AR(f)$ is free, see for instance \cite{DStFS}. 
Since the rank of the $S$-module $AR(f)$ is always $n$, it follows that $AR(f)$ is free exactly when $AR(f)$ is generated by $n$ homogeneous syzygies, say $\rho_1, \ldots, \rho_n$. We set  $d_j=\deg \rho_j$ and assume that 
$$ d_1\leq \ldots \leq d_n .$$
Note that the case $d_1=0$  occurs if and only if $D$ is a cone over a hypersurface in $\PP^{n-1}$, a situation which is not of interest in the sequel, hence we assume from now on $d_1 \geq 1$.
 We call these degrees $(d_1,\ldots,d_n)$ the {\it exponents} of the free hypersurface $D$ and  assume that $d_1 \leq \ldots \leq d_n$. Moreover, we  recall that 
 \begin{equation}
\label{eqF0}
 d_1+\ldots + d_n=d-1. 
  \end{equation}
  In particular, for $n=2$, one has
 \begin{equation}
\label{eq90}
d_1 \leq \frac{d-1}{2}.
\end{equation}
In the following result, we construct {\it explicit bases} for the $R$-module $Der(g)$ in many cases, in particular for the Milnor fibers of weighted homogeneous polynomials.
\begin{prop}
\label{propKS} 
Assume  that $h \in R$ is a weighted homogeneous polynomial of type  $(w_1,\ldots, w_n;N)$, where the weights $w_j$ are positive or negative integers, and the weighted degree $N$ is not zero.
Then the smooth affine hypersurface $X_t:g=h+t=0$, for $t \in \C^*$, is $A$-free. Moreover, when $w_1= \ldots=w_n=1$, the hypersurface $V:h=0$ in $\PP^{n-1}$ is free if and only if  the affine hypersurface $X$ is  $A_0$-free.
\end{prop}
Note that the hypersurfaces $X_t$ for $t \in \C^*$ are isomorphic to each other, and in particular $X_{-1}$ is usually called the {\it Milnor fiber of the weighted homogeneous polynomial } $h$.
\proof
For any non-zero polynomial $h \in R$, the exact sequence \eqref{es1} shows that the rank of the $R$-module $AR(h)$ is $n$. This implies that there are $n$ derivations $\be_j \in AR(h)$, for $j=1,\ldots,n$ such that the corresponding matrix $M(\be_1, \ldots, \be_n)$, constructed using their coefficients as in Proposition \ref{propA} above, has a non-zero determinant. In fact, we can and do choose $\be_1$ to be the Euler derivation $\epsilon$ as in the formula \eqref{Euler}, and $\be_j \in AR_0(h)$ for $j>1$. Consider now the following $n$ derivations $\delta_j \in AR(g)$, for $j=1,\ldots,n$, where
$$\delta_j=g_j\sum_{1\leq k\leq n, k \ne j}\frac{w_kx_k}{N}\partial_k+\left(t+\frac{w_jx_jg_j}{N}\right)\partial_j.$$
Note that one has $g_i=h_i$ for all $i=1,\ldots,n$ and
$$\delta_j g=\delta_jh=g_j(h+t)=g_jg.$$
Therefore we have indeed $\delta_j \in AR(g)$, for $j=1,\ldots,n$.
Consider now the matrix $M(\delta_1, \ldots, \delta_n)$ as in Proposition \ref{propA} and note that we have
$$M(\be_1, \ldots, \be_n) M(\delta_1, \ldots, \delta_n)=M(g\be_1, t\be_2,\ldots, t\be_n).$$
By taking determinants we get
$$\det M(\be_1, \ldots, \be_n) \det M(\delta_1, \ldots, \delta_n)= $$
$$=\det M(g\be_1,t\be_2, \ldots,t \be_n)=t^{n-1}g\det M(\be_1, \ldots, \be_n).$$
Since $\det M(\be_1, \ldots, \be_n) \ne 0$, we conclude that
$$ \det M(\delta_1, \ldots, \delta_n)=t^{n-1}g$$
and hence $AR(g)$ is a free $R$-module and $\delta_1, \ldots, \delta_n$ is a basis of it by Saito's Criterion in Proposition \ref{propA}.

To prove the claim about the $A_0$-freeness of $X_t$, note that $AR_0(g)=AR_0(h)$. Hence $X_t$ is $A_0$-free if and only if $V:h=0$ is a free hypersurface in $\PP^{n-1}$.
\endproof

The following result is an obvious consequence of Theorem \ref{thm1}.

\begin{cor}
\label{corF1}
Let $D:f=0$ be a  free hypersurface of degree $d>1$ in $\PP^n$. Consider an affine part of it, namely the affine hypersurface
$$ X_H=D \setminus H \subset \PP^n \setminus H \simeq \C^n,$$
for a hyperplane $H \subset \PP^n$, $H$ not an irreducible component of $D$. Then the affine hypersurface $X_H$ is  $A$-free.
\end{cor}
\proof
Since $D$ is free, the $S$-module $AR(f)$ admits a set of $n$ generators, and hence, by Theorem \ref{thm1}, the $R$-module $AR(g)$
also has a generating set of $n$ elements.
\endproof

The converse implication in Corollary \ref{corF1} does not hold, as the following example shows.

 \begin{ex}
\label{exF0} 
Consider the smooth projective surface $D: f=x_0^2+x_1^2+x_2^2+x_3^2=0$,
and the associated affine part $X:g=x_1^2+x_2^2+x_3^2+1$.
Then $D$ is not free, since a free projective hypersurface has a codimension 1 singular locus, see \cite{DStFS}. However, the affine surface $X$ is $A$-free by Proposition \ref{propKS}. More precisely
the $R$-module $AR(g)$ is generated by the vectors 
$$\rho_1'=(-2x_3,x_1x_3,x_2x_3,x_3^2+1), \  \rho_2'=(-2x_2,x_1x_2,x_2^2+1,x_2x_3) \text{ and } $$
$$ \rho_3'=(-2x_1,x_1^2+1,x_1x_2,x_1x_3).$$
On the other hand, the $R$-module $AR_0(g)$ is not free, since it is the first syzygy module of the maximal ideal $(x_1,x_2,x_3) \subset R$.
This example also shows that one cannot expect a formula similar to \eqref{eqF0} in the affine case. 
 \end{ex}
 \begin{rk}
\label{rkMac} 
 To check whether a given affine hypersurface of degree $d$ in $\C ^n$ is $A$-free or $A_0$-free, at least for small values of $d$ and $n \geq 3$, one may use {\it Macaulay2}, the package {\it QuillenSuslin}, see \cite{Mac}.
 \end{rk} 
For any $R$-module of finite type $M$, let $\widetilde M$ be the coherent sheaf on $\C^2$ associated to $M$. It is usual to say that an 
$R$-module of finite type $M$ has a property (P) if the associated
coherent sheaf $\widetilde M$ has this property.

The following result is at least folklore. Since we are not aware of any reference for it, we include a proof.

\begin{thm}
\label{thmA}
A reduced affine hypersurface $X$ is $A$-free if and only if all the hypersurface singularities $(X,p)$ for $p \in X_{sing}$ are free.
In particular, any smooth affine hypersurface $X$ is $A$-free.
\end{thm}
\proof

First we use Quillen-Suslin Theorem, saying that an $R$-module of finite type $P$ is free if and only if it is projective, see \cite{Q,Su}.
Then we use the general fact that an $R$-module of finite type $P$ is  projective if and only if all of its localizations $P_m$ at maximal ideals
$m \subset R$ are free, see \cite[Theorem A3.2]{Eis0}.

We apply these results to the $R$-module $Der(g)$. If the maximal ideal $m$ corresponds to a point $p_m \in \C^n$, then the germ $(X,p_m)$ of the hypersurface $X$ at $p_m$ can be regarded as a singularity at the origin of $\C^n$, simply by using a translation, so it has a local equation
$h^{p_m}=0$, with $h^{p_m} \in \OO_n$. When $p_m \notin X$, then one can take $h^{p_m}=1$, and hence
$$Der(h^{p_m})=Der(1)=Der(\OO_n)$$
is clearly a free $\OO_n$-module.
Similarly, when $p_m$ is a smooth point on $X$, then one can choose the local analytic coordinates such that $h^{p_m}=x_1$, and in this case
$Der(h^{p_m})=Der(x_1)$
is also a  free $\OO_n$-module, with a basis given by $$x_1\partial_1, \partial_2, \ldots, \partial_n.$$
For any maximal ideal $m \subset R$, we have a morphism of local rings
$$\iota_m: R_m \to \OO_n,$$
obtained by regarding a rational fraction in the localization $R_m$ as an analytic function germ at $p_m \in \C^n$, and then using the obvious translation to get an analytic function at the origin of $\C^n$, hence an element of $\OO_n$. Using $\iota_m$, one can regard the ring $\OO_n$ as an $R_m$-module.
Then one has 
$$Der(g)_m=Der(g) \otimes_R R_m \text{ and } Der(h^{p_m})=Der(g)_m \otimes_{R_m} \OO_n.$$
Now we can use general GAGA type results, see \cite{Serre}, and see that the $R_m$-module $Der(g)_m$ is free if and only if the $\OO_n$-module $Der(h^{p_m})$ is free. We have seen above that this latter condition is automatically satisfied when $p_m \notin X_{sing}$.
This completes the proof of our claims.
\endproof

\begin{cor}
\label{corA} 
Any reduced affine plane curve $X$ is both $A$-free and $A_0$-free.
Moreover, one has either $A(X)=R$, or the proper ideal $A(X)\subset R$ is minimally generated by exactly two polynomials.
\end{cor} 

\proof
The fact that $Der_0(X)$ is a free $R$-module of rank 1 follows from 
Example \ref{exC1}.  For $Der(X)$ we use Theorem \ref{thmA} and hence it is enough to show that any plane curve singularity is free. To do this, one can use two approaches.
Namely, one can use K. Saito results, see in particular  \cite[Corollary 1.7]{KS}. Or, one can notice that the analytic sheaf version of the exact sequence \eqref{es1} implies that  $Der(h)$ is a stalk of a second syzygy sheaf $\F$, see for instance \cite[Definition 1.1.5]{OSS}. Hence the codimension of the singularity set $S(\F)$ of $\F$ is $>2$, see \cite[Theorem 1.1.6]{OSS}. Therefore in our case $S(\F)=\emptyset$ and the coherent sheaf $\F$ is locally free. Since $Der(h)=\F_0$, the stalk of $\F$ at the origin, then the first claim is proved.

Theorem \ref{thmA} and the exact sequence \eqref{ideal1.2} imply that the ideal $A(X)$ can be generated by 2 elements. Since $A(X)$ has finite codimension in $R$, it follows that $A(X)$ is generated by a single element if and only if this element is a unit, i.e. exactly when $A(X)=R$.
\endproof

 \begin{ex}
\label{ex0} 
Let $C:f=0$ be a smooth curve of degree $d>1$ in $\PP^2$. Then the projective curve $C$ is not free, see for instance \cite[Theorem 2.8]{DStFD}, but all the affine parts of it, namely the curves
$$ X_L=C \setminus L \subset \PP^2 \setminus L \simeq \C^2,$$
for all lines $L \subset \PP^2$ are free in view of Theorem \ref{thmA}. 
For instance, if $C:f=x_0^d+x_1^d+x_2^d=0$ and $X:g=x_1^d+x_2^d+1=0$, then $AR(g)$ is generated by
$$\rho'_1=(-dx_2^{d-1}, x_1x_2^{d-1},x_2^d+1) \text{ and } \rho'_2=(-dx_1^{d-1},x_1^d+1,x_1^{d-1}x_2).$$
To check this claim, we recall that for a smooth curve $C:f=0$, the $S$-module $AR(f)$ is generated by the Koszul syzygies $\kappa_{ij}$ defined in \eqref{eq7.05}. It follows from Theorem \ref{thm1} that the $R$-module $AR(g)$ is generated by
the images $\kappa'_{ij}$  of the Koszul syzygies under the morphism $\phi$. 
The reader can check easily that  one has
$\kappa'_{02} =d\rho_1'$, $\kappa'_{01}= d\rho_2'$ and 
$\kappa'_{12}=dx_1^{d-1}\rho_1' -dx_2^{d-1}\rho'_2.$ 
Note that $sdeg(\rho'_1)=sdeg(\rho'_2)=d-1$, and in particular the minimal degree $mdr(g)$ of a generator for $AR(g)$ for a free affine curve does not satisfy an inequality similar to \eqref{eq90}. For a general smooth curve $C$ and any line $L$, the equality $mdr(g)=mdr(f)=d-1$ follows from Corollary \ref{corthm1}.
\end{ex}

\begin{rk}
\label{rk0.2}
For any pair of curves $X:g=0$, $C:f=0$ as above, let $r=mdr(f)=mdr(g)$. By Corollary \ref{corthm1} we know that $\dim AR(f)_r=\dim AR(g)_{\leq r}$.
Note that the integer $\dim AR(f)_r$ is the number of generators of degree $r$ in any minimal system of homogeneous generators for the graded $S$-module $AR(f)$. In particular, $\dim AR(f)_r \leq N_f$, where $N_f$ is the minimal number of generators for the graded $S$-module $AR(f)$. The equality can hold, for instance when $C:f=0$ is a maximal nodal curve as in \cite{Oka}, namely a rational curve of degree $d$ having exactly $(d-1)(d-2)/2$ nodes. In this case $r=d-1$ and 
$$\dim AR(f)_r =N_f=d+1,$$
see \cite[Proposition 5.8]{DStCol}. In particular, the difference
$N_f-N_g=d-1$ can be as large as we wish.
On the other hand, the inequality $\dim AR(g)_{\leq r} \leq N_g$, where $N_g=2$ is the minimal number of generators for the  $R$-module $AR(g)$ can fail. For instance, for a maximal nodal curve as above, we have $\dim AR(g)_{\leq r}=d+1$, but $N_g=2$.

\end{rk}

\section{Bounds on Tjurina numbers} 

In this section we consider affine hypersurfaces $X$ having only isolated singularities. We discuss first the case of plane curves, where the results are more complete.
Let $\tau(C)$ denote the total Tjurina number of a projective plane curve $C:f=0$, which is by definition the sum of all local Tjurina numbers $\tau(C,p)$ for $p \in C_{sing}$. The two invariants $\tau(C)$ and $mdr(f)$ are related by the following inequalites, see \cite{Dmax,dPW,E}.
\begin{thm}
\label{thmtau2} 
Let $C:f=0$ be a reduced projective plane curve of degree $d$ with $r=mdr(f)$.
Then one has
$$\tau(d,r)_{min}= (d-1)(d-r-1) \leq \tau(C).$$
On the other hand, one has the following upper bounds for $\tau(C)$.
\begin{enumerate}

\item If $r <d/2$, then $$\tau(C) \leq \tau(d,r)_{max}= (d-1)(d-r-1)+r^2$$ and the equality holds if and only if the curve $C$ is free.

\item If $d/2 \leq r \leq d-1$, then
$$\tau(C) \leq \tau(d,r)_{max},$$
where, in this case, we set
$$\tau(d,r)_{max}=(d-1)(d-r-1)+r^2-{ 2r-d+2 \choose 2}.$$

\end{enumerate}

\end{thm} 
The equality holds  in (2) if and only if the projective curve $C$ is maximal Tjurina, see  \cite{DStCol} for the definition and the properties of these curves. 
In the affine case, we have the following result, similar to Theorem \ref{thmtau2}.
\begin{cor}
\label{corthm2}
For any pair of curves $X:g=0$, $C:f=0$ of degree $d$ as above, one has
$$\tau(X) \leq \tau(C)$$
and the equality holds exactly when the curve $C$ has no singularities on the line at infinity $L_0:x_0=0$.
Moreover  
$$\tau(X) \leq \tau(d,r)_{max},$$
where $r=mdr(g)=mdr(f)$ and $\tau(d,r)_{max}$ is defined as in Theorem \ref{thmtau2}.
In addition, this last inequality is an equality if and only if the curve $C$  has no singularities on the line at infinity $L_0$ and $C$ is either  free when $r'<d/2$, or maximal Tjurina when $r'>d/2$.
\end{cor}
\proof
The first claim is obvious. To prove the remaining claims, we use Theorem \eqref{thmtau2} and the equality $mdr(f)=mdr(g)$ established in Corollary \ref{corthm1}.

 \endproof

In higher dimensions, one has the following result, see \cite[Theorem 5.3]{duPCTC01} for the proof and \cite{DAndr} for a discussion of this result.
\begin{thm}
\label{thmC}
 Assume that the hypersurface $D:f=0$ in $\PP^n$ has only isolated singularities and set $r=mdr(f)$.
Then
$$(d-r-1)(d-1)^{n-1} \leq \tau(D) \leq \tau(d,r)_{max}=(d-1)^n-r(d-r-1)(d-1)^{n-2}.$$
\end{thm}
Using this result, we have the following analog of Corollary \ref{corthm2}.

\begin{cor}
\label{corthmC}
For any pair of hypersurfaces  $X:g=0$ in $\C^n$ and  $D:f=0$ in $\PP^n$ of degree $d$ as above, where $D$ has only isolated singularities, one has
$$\tau(X) \leq \tau(D)$$
and the equality holds exactly when the hypersurface $D$ has no singularities on the hyperplane at infinity $H_0:x_0=0$.
Moreover  
$$\tau(X) \leq \tau(d,r)_{max},$$
where $r=mdr(g)=mdr(f)$ and $\tau(d,r)_{max}$ is defined as in Theorem \ref{thmC}.
\end{cor}

One may ask if the degree $mdr(g)$ (resp. $mdr(f)$) depends on the singularities of $X$ (resp. $D$). One has the following partial result, when $D$ has only isolated singularities, saying that the answer is negative  when $\tau(X)$ (resp. $\tau(D)$) is small compared to $(d-1)^{n-1}$.

\begin{prop}
\label{propB} 
For any pair of hypersurfaces $X:g=0$ in $\C^n$ and  $D:f=0$ in $\PP^n$ of degree $d$ as above,  where $D$ has only isolated singularities, one has
$$d-1-\frac{\tau(D)}{(d-1)^{n-1}} \leq mdr(g)=mdr(f) \leq d-1.$$
In particular, if in addition $\tau(D) <(d-1)^{n-1}$, then $mdr(g)=mdr(f)=d-1$.
\end{prop}
\proof
The first inequality follows from the lower bound for $\tau(D)$ in Theorem \ref{thmC}. The relations $mdr(g)=mdr(f) \leq d-1 $ come from 
Corollary \ref{corthm1}. 
\endproof
Note that in Proposition \ref{propB} one can replace $\tau(D)$ by $\tau(X)$, when the affine hypersurface $X$ has no singularities at infinity, namely when the hypersurface $D$ has no singularities on the hyperplane at infinity $H_0:x_0=0$.

\end{document}